\documentclass[11pt]{amsart}
\usepackage{cite}

\usepackage{bm}
\usepackage[all]{xy}
\usepackage[utf8]{inputenc}

\usepackage[english]{babel}
\usepackage{stmaryrd}

\newtheorem{thmx}{Theorem}
 
\newtheorem{theorem}{Theorem}[section]
\newtheorem{definition}[theorem]{Definition}
\newtheorem{lemma}[theorem]{Lemma}
\newtheorem{corollary}[theorem]{Corollary}
\newtheorem{remark}[theorem]{Remark}
\newtheorem{proposition}[theorem]{Proposition}

\newcommand{\al}{\alpha}

\newcommand{\de}{\delta}
\newcommand{\e}{\varepsilon}

\newcommand{\Lm}{\Lambda}

\newcommand{\p}[1]{\partial^{[#1]}}
\newcommand{\HH}{{\mathbb{H}}}

\newcommand{\CC}{{\mathbb{C}}}
\newcommand{\PP}{{\mathbb{P}}}

\newcommand{\ZZ}{{\mathbb{Z}}}
\newcommand{\RR}{{\mathbb{R}}}
\newcommand{\op}{\operatorname}

\newcommand{\T}{\Theta}

\newcommand{\ch}[1]{\left[\begin{smallmatrix}#1\end{smallmatrix}\right]}
\newcommand{\tch}[1]{{\vartheta\left[\begin{smallmatrix}#1\end{smallmatrix}\right]}}
\newcommand{\sm}[1]{\left(\begin{smallmatrix}#1\end{smallmatrix}\right)}
\newcommand{\kum}[1]{K_{#1}}

\newcommand{\mc}{\mathcal}

\def\Oo{\mathcal O}

\def\C{\mathbb C}

\def\p{\mathbb P}
\def\cc{\mathcal C}

\def\cc{\mathcal{C}}

\newcommand{\m}{\mbox}

\newcommand{\xx}{\otimes}

\newcommand{\beginCD}{\begin{equation*}\begin{CD}}
\newcommand{\enCD}{\end{CD}\end{equation*}}

\DeclareMathOperator{\Aut}{Aut}

\DeclareMathOperator{\id}{id}

\DeclareMathOperator{\sym}{Sym}

\DeclareMathOperator{\Ho}{Hom}

\DeclareMathOperator{\Spec}{Spec}
\DeclareMathOperator{\ext}{Ext}

\DeclareMathOperator{\Pic}{Pic}

\DeclareMathOperator{\HOM}{HOM}

\newcommand{\K}{\mathcal K}
\newcommand{\A}{\mathcal A}
\DeclareMathOperator{\abs}{abs}
\DeclareMathOperator{\pre}{pre}
\DeclareMathOperator{\emb}{em}
\newcommand{\Ka}{\K^{\abs}}
\newcommand{\Ke}{\K^{\emb}}
\def\kk{{\bm{\mathsf{K}}}^{\emb}}
\def\aa{{\bm{\mathsf{A}}}}

\def\H{\mathbb H}
\def\hh{\mathcal H}
\DeclareMathOperator{\hilb}{Hilb}
\DeclareMathOperator{\isom}{Isom}
\newcommand{\Isom}{\pmb{\isom}}
\newcommand{\Hoo}{\pmb{\Ho}}
\newcommand{\zzz}{\underline{\ZZ\slash 2\ZZ}}
\newcommand{\Ag}{\A_g \negmedspace\negmedspace\fatslash(\zzz)}

\begin{document}

\markboth{Mattia Galeotti, Sara Perna}{Moduli spaces of abstract and embedded Kummer varieties}

\title{Moduli spaces of abstract and embedded Kummer varieties}

\author{Mattia Galeotti}

\author{Sara Perna}

\maketitle

\begin{abstract}
In this paper, we investigate the construction of two moduli stacks of Kummer varieties. The first one is the stack $\mathcal K^{\text{abs}}_g$ 
of abstract Kummer varieties and the second one is the stack $\mathcal K^{\text{em}}_g$ of embedded Kummer varieties. 
We will prove that $\mathcal K^{\text{abs}}_g$ is a Deligne-Mumford stack and its coarse moduli space is isomorphic to $\boldsymbol A_g$, 
the coarse moduli space of principally polarized abelian varieties of dimension $g$. On the other hand we give a modular family $\mathcal W_g\to U$ 
of embedded Kummer varieties embedded in $\mathbb P^{2^g-1}\times\mathbb P^{2^g-1}$, meaning that every geometric fiber of this family 
is an embedded Kummer variety and every isomorphic class of such varieties appears at least once as the class of a fiber. 
As a consequence, we construct the coarse moduli space $\boldsymbol{\mathsf K}^{\text{em}}_2$ of embedded Kummer 
surfaces and prove that it is obtained from $\boldsymbol A_2$ by contracting the locus swept by a particular linear equivalence class of curves. 
We conjecture that this is a general fact: $\boldsymbol{\mathsf K}^{\text{em}}_g$ could be obtained from $\boldsymbol A_g$ via a contraction for all $g>1$.
\end{abstract}


\section{Introduction}
The main subject of this work are the moduli of Kummer varieties.
These varieties can be defined in an abstract way as the quotient 
of an abelian variety, or as the image of a natural
embedding from principally polarized abelian varietis.

In order to define the moduli of abstract Kummers,
let $A$ an abelian variety of dimension $g$ and denote by $i:A\to A$ its natural involution. We define the abstract Kummer variety of $A$ as $\kum{A}:= A/\langle i\rangle$. If $g\ge2$, $\kum A$ is a singular variety of the same dimension as $A$ with $2^{2g}$ singular points 
corresponding to the image of the 2-torsion points via the projection $\pi:A\to A/\langle i\rangle$ (see~\cite[Section 4.8]{BL04}). 

The first result of our paper concerns the moduli stack $\Ka_g$ of
abstract Kummer varieties of dimension $g$, that is the stackification of the category
fibered in groupoids over schemes whose objects are 
abelian schemes quotiented by their natural involution (see Definition \ref{def_ka}).

\begin{thmx}
The category $\Ka_g$ is a Deligne-Mumford stack. Its
coarse space is the moduli space $\aa_g$ of principally polarized
abelian varieties.
\end{thmx}

We consider now an abelian variety over $\C$, carrying a principal polarization~$\T$.
Note that the line bundle $2\T$ is canonically defined. Indeed, one can choose a symmetric divisor $D$ in the class of $\T$ and then the line bundle $2D$ will not depend on $D$. 
Let $\varphi:A\to\PP^{2^g-1}$ be the map defined by the line bundle $2\T$. 
We recall that, with a suitable choice of a theta structure, the map $\varphi$ is given by second order theta functions. 
If the principally polarized abelian variety $(A,\Theta)$ is indecomposable, then $\varphi$ factors through the quotient $\pi: A\to \kum A$ followed by an embedding of $\kum{A}$ into $\PP^{2^g-1}$.
If $(A,\Theta)$ is decomposable and $(A_1,\Theta_1)\times\dots\times (A_s,\Theta_s)$ is the decomposition of $(A,\Theta)$ into irreducible principally polarized abelian varieties of lower dimension, then 
$\varphi$ factors through the product 
$\kum{A_1}\times\dots\times\kum{A_s}$ and moreover $\varphi$
 gives an embedding of this product into $\p^{2^g-1}$ (see~\cite[Section 4.8]{BL04}).
The image of $\varphi$ is a $g$-dimensional variety of degree $(2\T)^g/2^s=2^{g-s}g!$. 
Clearly, if a principally polarized abelian variety is indecomposable, then its associated abstract  Kummer variety is isomorphic
to the image of the corresponding $\varphi$ morphism.\newline

In the second part of our work we focus on the case of embedded Kummer surfaces over~$\C$.
In order to define the moduli of embedded Kummers, in \S\ref{sec_emb} we introduce 
a family 
$$\xymatrix{\mc W_2\ar[d] \ar@{^{(}->}[r]& U\times \p^{3}\\
U\\}$$
defined via theta functions on the Siegel space $\HH_2$. 
An embedded Kummer surface is therefore defined as a geometric fiber of this family.
The moduli stack $\Ke_2$ of embedded Kummer surfaces is the stackification
of the category fibered in groupoids over schemes whose objects 
are flat families Kummer surfaces embedded 
in a projective bundle (see Definition~\ref{def_ke}).

The abstract Kummer surface of an indecomposable complex abelian surface over is a quartic hypersurface in $\PP^3$ with 16 nodes.
These have been completely classified in terms of configurations of points and planes in $\PP^3$ (see~\cite{GD}). The blow up of this Kummer surface at the 16 nodes is known to be a $K3$ surface.
In the case of a decomposable abelian surface $E_1\times E_2$, where $E_1$ and $E_2$ are two elliptic curves, the associated Kummer variety
is a double quartic. We denote by $\aa_{1,1}$ the image of the natural map $\aa_1\times \aa_1\to \aa_2$
sending a pair of elliptic curves to their product.
If we fix any elliptic curve $B$, we denote by $C_1^{(B)}$ the curve in $\aa_2$ which is the locus
of abelian surfaces $B\times E$ where $E$ is any elliptic curve. The linear equivalence class of $C_1^{(B)}$ is independent of $B$
and is denoted by~$C_1$. We denote by $\aa_2\to \aa_2^0$ the Mori contraction on $\aa_2$
associated with the ray $\RR_+C_1$, this means that  the entire locus $\aa_{1,1}$ gets contracted. 

\begin{thmx}
The category $\Ke_2$ of embedded Kummer surfaces  is a Deligne-Mumford stack.
Its associated moduli space  $\kk_2$ is isomorphic to the space $\aa_2^0$, 
obtained via the Mori contraction above.
\end{thmx}

 We conjecture that this is a general fact, that $\kk_g$ should be
obtained from $\aa_g$ via a contraction  for all $g>1$. A similar statement has
been proved in \cite{SB16} for the perfect compactification $\aa_g^P$
and $g\leq 11$.\newline

Let us outline the structure of the paper.
In Section \ref{ppav}, we introduce principally polarized abelian varieties,
the moduli stack $\A_g$ and the Siegel space

 $\HH_g$.  In Section \ref{abs} ,
 we treat functor $\Ka_g$ of abstract Kummer varieties,
 while in Section \ref{emb} we introduce the universal embedded Kummer surface.
In Section \ref{embsurf}, we introduce the embedded Kummer functor $\Ke_2$ and
construct its coarse moduli space $\kk_2$.
Finally, in Section \ref{degen},
we consider rank $1$ degenerations of principally polarized abelian varieties aiming at pointing to a possible (partial) compactification of the moduli space of embedded Kummer varieties.

\section{Preliminaries}\label{ppav}

In this section, we  briefly give the definition of principally polarized abelian varieties over an algebraically
closed field $k$ and describe their moduli space, then we give the general setup more specifically for the complex case.

\subsection{Principally polarized abelian varieties}
An abelian scheme over a scheme $S$ is an $S$-group scheme $A\to S$ which is
proper, flat, finitely presented and has smooth and connected geometric fibers.
A well known fact is that such an abelian scheme is actually a commutative group scheme.
As a consequence there exists a $0$-section $\varepsilon\colon S\to A$. For any invertible sheaf $\mc L$
over $A$, a rigidification is an isomorphism $\Oo_S\to \varepsilon^*\mc L$.

For every abelian scheme $A\to S$ we define its dual abelian scheme $\hat A\to S$.
Consider the relative Picard functor $\pmb\Pic(A\slash S)$, which is the functor
from $Sch\slash S$ to sets such that
$$\pmb\Pic(A\slash S)(T)=\left\{\begin{array}{c}
     \m{isomorphism classes of line bundles }\mc L\m{ on }\\
     A_T=A\times_S T\m{ rigidified along }\varepsilon_T=\varepsilon\times_S T
     \end{array}
\right\}.$$
The Picard functor is represented by a scheme
locally of finite presentation over $S$ denoted by $\Pic(A\slash S)$.
We denote by $\pmb\Pic^0(A\slash S)$ the subfunctor of $\pmb\Pic(A\slash S)$
corresponding to line bundles $\mc L$ such that for all $t\in T$, $\mc L\xx k(t)$ is
algebraically equivalent to zero on $A_t$. This functor is represented
by a subscheme $\Pic^0(A\slash S)$ of $\Pic(A\slash S)$ which
is an abelian scheme over S (see \cite[\S1]{FC}), called the dual abelian scheme
of $A$ and denoted by $\hat A:=\Pic^0(A\slash S)$. There exists
a canonical isomorphism $A=\hat{\hat{A}}$.

If $\pi\colon A\to S$ is the structural morphism of $A$, any invertible sheaf $\mc L$
on $A$ defines a morphism
$\phi_{\mc L}\colon A\to \hat A$
such that for all $a\in A(S)$,
$$a\mapsto T_a^*\mc L\xx\mc L^{-1}\xx (a\circ \pi)^*\mc L^{-1}\xx(\varepsilon\circ \pi)^*\mc L,$$
where $T_a\colon A\to A$ is the translation by $a$.
The morphism $\phi_{\mc L}$ is in fact a group morphism
by the theorem of the cube (see \cite[Chap.3, p.91]{Mum70}). If $\mc P$ is the
Poincaré line bundle over $A\times_S \hat A$ rigidified
along  $\hat A\cong S\times_S\hat A\xrightarrow{\varepsilon\times \id_{\hat A}}A\times_S \hat A$
and $\upsilon\colon A\to \hat A$
is an abelian scheme morphism, we define a new line bundle over~$A$,
$\mc L^\Delta(\upsilon):=(\id_A\times \upsilon)^*\mc P$.
Furthermore, the following is a well-known fact (see again \cite[\S I.1]{FC})
\begin{equation}\label{eq_lambda}
\phi_{\mc L^\Delta(\upsilon)}=2\upsilon.\newline
\end{equation}

If $S=\Spec(k)$ is the spectrum of a field, an abelian scheme over $S$
is called an abelian variety over $k$.
\begin{definition}
Given an abelian scheme $A\to S$, a polarization is an abelian scheme morphism
$\upsilon\colon A\to \hat A$ such that for each geometric point $\bar s$ of $S$
there exists an ample line bundle $L_{\bar s}$ over $A_{\bar s}$ with the property that
$\upsilon_{\bar s}=\phi_{L_{\bar s}}$.
The polarization is principal if it is an isomorphism.
\end{definition}

Note that the ample line bundle $L_{\bar s}$, on each geometric fiber, appearing
in the above definition is unique only up to algebraic equivalence.

We introduce now the moduli problem of principally polarized abelian schemes, using
the language of categories fibered in groupoids over schemes.

\begin{definition}
Let $g\geq 1$. We denote by $\A_g$ the category fibered in groupoids whose objects are principally
polarized abelian schemes $(A\to S, \upsilon)$ of relative dimension $g$, and whose
morphisms are abelian scheme isomorphisms commuting with the polarizations.
The fibration functor $\A_g\to Sch$ is the natural forgetful map. For any scheme $S$, the fiber $\A_g(S)$
is the groupoid of abelian schemes over the scheme $S$ and whose maps are lifting of $\id_S$.
\end{definition}

The category $\A_g$ is in fact a stack, and more notably a Deligne-Mumford stack: see for example \cite{FC}
for a proof of this.

\subsection{Complex tori}

An abelian variety over $\CC$ is a complex torus. Moreover, since
 the elements of the Néron-Severi group $\op{NS}(A)$ are first Chern classes of line bundles on $A$,
a polarization on $A$ is given by a positive definite Hermitian form $H\in\op{NS}(A)$.
We will denote by $(A,L)$ or $(A,H)$ the abelian variety $A$ with polarization given by (the first Chern class of) a line bundle $L$ or with polarization given by a positive definite $H\in\op{NS}(A)$, respectively.

Let $A=V/\Lm$ be a complex torus. 
We denote by $\Pi$ the period matrix of~$A$ 
\begin{theorem}[Riemann's relations,  see~{\cite[Theorem 4.2.1]{BL04}}]
A complex torus $A$ with period matrix $\Pi$ is an abelian variety if and only if there is a non degenerate alternating matrix $M\in M_{2g}(\ZZ)$ such that
\[\Pi M^{-1} {}^t\Pi=0\ \ \text{and}\ \  i\Pi M^{-1} {}^t\overline{\Pi}>0.\]
\end{theorem}

Suppose $A$ is an abelian variety with polarization of type $D$. 
With respect to opportune symplectic bases, the period matrix of $A$ is of the form $\Pi=(\tau,D)$ for some $\tau\in M_g(\CC)$
and the matrix $\tau$ 
is in the Siegel space of degree $g$, 
\begin{equation}\label{Siegel_space}
\HH_g=\{\tau\in M_g(\CC)\mid {}^t\tau=\tau,\,\op{Im}(\tau)>0\}.
\end{equation}
If one fixes the type $D$ of the polarization, any abelian variety with a polarization of type $D$ together 
with the choice of a symplectic basis of its lattice defines an element of $\HH_g$. 
Conversely, any matrix $\tau$ in the Siegel space defines a polarized abelian variety of type $D$ with a symplectic 
basis, $A_\tau=\CC^g/\tau\ZZ^g\oplus D\ZZ^g$.

Let $\Gamma_g$ be the integral symplectic group of degree $g$:
\[\Gamma_g=\left\{\gamma\in M_{2g\times 2g}(\ZZ)\mid \gamma\sm{0&\mathbf{1}_g\\-\mathbf{1}_g&0}{}^t\gamma=\sm{0&\mathbf{1}_g\\-\mathbf{1}_g&0}\right\}.\]
It acts on $\HH_g$ by a generalization of the M\"obius transformations of the complex plane. If $\gamma=\sm{a&b\\c&d}\in \Gamma_g$, where $a,\,b,\,c,\,d$ are $g\times g$ matrices, then for any $\tau\in\HH_g$ the action is defined as
\[\gamma\cdot\tau:=(a\,\tau+b)(c\tau+d)^{-1}.\]
It is classically known that for $\tau,\,\tau'\in\HH_g$ the principally polarized abelian varieties $A_\tau$ and $A_{\tau'}$ are isomorphic if and only if there exists $\gamma\in\Gamma_g$ such that $\tau'=\gamma\cdot\tau$. So the points of 
\begin{equation}\label{moduli_ppav}
\aa_g:=\HH_g\slash\Gamma_g
\end{equation}
are in one to one correspondence with the set of principally polarized abelian varieties of dimension $g$. Indeed, $\aa_g$ is the coarse moduli space of these varieties (see \cite[\S1.6]{FC}).

\section{The abstract Kummer stack}\label{abs}

The first Kummer stack we introduce is the stack of abstract Kummer varieties
of dimension~$g$,
denoted by $\Ka_g$. We will define the category fibered in groupoids
of families of abstract Kummer varieties, and $\Ka_g$ will be its stackification.
In this section we work over a field $k$ of characteristic different from $2$.

\begin{definition}\label{def_ka}
For any $g\geq 2$, $\K^{\abs,\pre}_g$ is a category fibered in groupoids over $Sch\slash k$ whose
objects are 
morphisms $K\to S$
such that there exists an abelian scheme $A\to S$ of relative dimension $g$,
with $K=A\slash i$.

We will call any pair  $K\to S$ as above a Kummer scheme over $S$. 
Given any abelian scheme $A\to S$, the Kummer scheme $K=A\slash i\to S$ 
defined as above will be called the Kummer scheme associated to $A\to S$.

A morphism in $\K^{\abs, \pre}_g$ between two Kummer schemes 
$K'\to S'$ and $K''\to S''$
is a pair of morphisms $(u\colon K'\to K''; v\colon S'\to S'')$,
such that the following diagram is a cartesian square
$$
\xymatrix{
K'\ar[r]^u \ar[d]& K''\ar[d]\\
S'\ar[r]^v & S''}
$$
and, moreover, such that there exists two abelian schemes $A', A''$
such that $K'= A'\slash i'$, $K''=A''\slash i''$ and $u$ lifts to a
scheme morphism $A'\to A''$. 
\end{definition}

We work in $g\geq2$, because with this constraint
the set of $K\to S$ ramification points coincides with the set
of singular points of $K$.
In the case $g=1$ in order to properly define $\K_1^{\abs,\pre}$
we have to introduce a section $\sigma\colon S\to \sym_S^{4}(K)$
sending any point to the ramification set of its fiber. With
this correction the results of this section remain true
for $g=1$ too.\newline

We will prove later that the category $\K_g^{\abs,\pre}$ is
in fact a prestack, and this justifies the notation we choose.
Before that, we will only consider $\K_g^{\abs,\pre}$ as a 
category fibered in groupoids over $Sch\slash k$.

\begin{definition}
The stack $\Ka_g$ of abstract Kummer varieties
of dimension $g$, is the stack associated to $\K_g^{\abs,\pre}$, i.e.~it is a stack equipped
with a morphism of categories $F\colon \K_g^{\abs,\pre}\to\Ka_g$ and such that for any stack $\K'$
the functor
$$\HOM(\Ka_g,\K')\to \HOM(\K_g^{\abs,\pre},\K')$$
induced by $F$, is an equivalence of categories.
\end{definition}

Given any category fibered in groupoids $\cc$ over $Sch\slash k$, two objects $X'$ and $X''$ over
the schemes $S'$ and $S''$ respectively and a morphism $f\colon S \to S'$, we denote by $\isom_{\cc,f}(X',X'')$ the
set of morphisms $X'\to X''$ over $f$.

Given an abelian scheme $A\to S$, we define the abstract group
$$\Aut_S(A)=\isom_{Sch\slash S, \id_S}(A,A).$$
Furthermore we define the group scheme $\zzz$ such that
if $S$ is any scheme with $V(S)$ as its set of connected components,
then $\zzz(S):=(\ZZ\slash 2\ZZ)^{\times V(S)}$. For every abelian scheme $A\to S$
there exists a natural inclusion
$$\iota_{A\slash S}\colon \zzz(S)\hookrightarrow \Aut_S(A).$$
Indeed, it suffices to define $\iota_{A\slash S}$ in the case of $S$ a connected scheme,
and in this case $\iota_{A\slash S}$ sends the non-trivial element of $\zzz(S)$ to
the canonical involution in $\Aut_S(A)$.

The inclusion $\iota_{A\slash S}$ is compatible with pullbacks
in the sense of \cite[p.16]{ACV03}. Indeed, 
for every base change $f\colon T\to S$ and any morphism of abelian schemes $\phi\colon f^*A\to A$,
we have
$$\iota_{f^*A\slash T}\circ f^*=\phi^*\circ \iota_{A\slash S}\colon \zzz(S)\to \Aut(f^*A).$$

As a consequence
we can introduce the rigidification $\Ag$, a stack with some important properties.

\begin{theorem}[See {\cite[Theorem 5.1.5]{ACV03}}]\label{teo_this}
There exists a stack $\Ag$ and a natural morphism $F\colon \A_g\to\Ag$ 
with the following properties:
\begin{enumerate}
\item $F$ is a $\zzz$-gerbe;
\item  for
all objects $A\to S$ in $\A_g(S)$ any element of $\zzz(S)\subset \Aut_S(A)$ maps to the identity via $F$, 
and $F$ is universal with
respect to this property, i.e.~any stack morphism with this same property
factors through $F$;
\item $\Ag$ is Deligne-Mumford;
\item $\Ag$ has a coarse moduli space which is isomorphic to $\aa_g$.
\end{enumerate}
\end{theorem}

As stated in the proof of the paper of Abramovich-Corti-Vistoli \cite{ACV03},
the stack $\Ag$ is the stackification of the prestack $(\Ag)^{\pre}$,
who has the same objects of~$\A_g$. Moreover, if we consider
two abelian schemes $A'\to S', A''\to S''$, then
for any map $f\colon S'\to S''$ the set of morphisms $A'\to A''$ over $f$ is
\begin{equation}\label{eq_agpre}
\isom_{(\Ag)^{\pre},f}(A',A'')=\frac{\isom_{\A_g,f}(A',A'')}{\zzz(S'')},
\end{equation}
where $\zzz(S'')\subset\Aut_{S''}(A'')$ acts naturally on $\isom_{\A_g,f}(A',A'')$.
This justifies the notion of rigidification, because $\Ag$ is obtained
by keeping the same objects -- \textit{i.e.}~abelian schemes --
and quotienting out the natural abelian scheme involution, over the whole $\A_g$ stack.\newline

Given any abelian scheme $A\to S$, 
and two elements $s_1,s_2\colon S\to A$, we denote by
$\mu(s_1,s_2)=s_1+s_2$ the multiplication morphism.
We consider the canonical $0$-section $\varepsilon\colon S\to A$,
multiplication $\mu\colon A\times_SA\to A$, the diagonal
morphism $\Delta\colon A\to A\times_S A$ and
the involution $i\colon A\to A$.
We define the subscheme of $2$-torsion points
$$A[2]:=(\mu\circ \Delta)^{-1}\left(\varepsilon(S)\right),$$
and we denote by $\Ho_S(S,A[2])$ the set of sections of $A[2]\to S$.\newline

Consider the morphism of categories
fibered in groupoids over schemes
$$G^{\pre}\colon (\Ag)^{\pre}\to\K_g^{\abs,\pre},$$
that sends any abelian scheme $A\to S$ to its associate Kummer scheme $K=A\slash i\to S$.
By the definition of $\K_g^{\abs,\pre}$, the morphism $G^{\pre}$ is essentially surjective.
Moreover, by the universal property of stackification, it extends to a morphism of stacks
$G\colon\Ag\to\K_g^{\abs}$ that fits into the following commutative diagram.
\begin{equation}\label{diagG}
\xymatrix{
\A_g\ar@{->>}[r]\ar@{->>}[d] &\K_g^{\abs,\pre}\ar[r] &\K_g^{\abs}\\
(\Ag)^{\pre}\ar@{->>}[ru]_{G^{\pre}}\ar[r]& \Ag\ar[ru]_{G}\\
}
\end{equation}
The morphisms marked with a double arrow are essentially surjective, and the two
other horizontal arrows are the stackification morphisms.\newline

Given any category $\cc$ fibered in groupoids over the category of schemes, if we have two objects $x,y$ over
the same base scheme $S$, then we denote by 
$$\isom_\cc(x,y):=\isom_{\cc,\id_S}(x,y)$$ 
the set
of isomorphisms over the identity morphism~$\id_S$. If the category $\cc$
is clear from context, we will simply write $\isom(x,y)$. Furthermore, we introduce
the functor $\Isom_\cc(x,y)$: this is a functor from $Sch\slash S$ to sets that sends
any scheme morphism $w\colon T\to S$ to the set $\isom_T(w^*x,w^*y)$.\newline

Recall that, given any abelian scheme $A$, $A[2]$ is the
group subscheme of $2$-torsion points. We observe that 
the group scheme $A[2]$ is also well defined for any object $A$
of $\Ag$.
Observe furthermore, that the translation $t_a\colon A\to A$ by a $2$-torsion point $a\in A[2]$,
pushforwards to the associated abstract Kummer $K_A$.

Consider two abelian schemes $A'\to S$ and $A''\to S$,
and the associated
Kummer schemes $K'\to S$ and $K''\to S$.
The next lemma gives a relation between $\Isom_{(\Ag)^{\pre}}(A',A'')$ and
$\Isom_{\K_g^{\abs,\pre}}(K',K'')$, where both sets are intended over the same base
scheme $S$, but in different categories.

\begin{lemma}\label{lemma_aut}	
We have the equality of functors
$$\Isom_{\K_g^{\abs,\pre}}(K',K'')=\Hoo_S(-,A''[2])\times\Isom_{(\Ag)^{\pre}}(A',A''),$$
where $\Hoo_S(-,A''[2])$ is the functor sending any object $T$ of $Sch\slash S$ to the
set of morphisms $\Ho_S(T,A''[2])$.
\end{lemma}

The above lemma means that an isomorphism of $K'$ and $K''$
is the composition of a $2$-torsion point translation, with
an isomorphism $A'\to A''$ in $\Ag^{\pre}$. Equivalently,
Isomorphisms $K'\to K''$ up to $2$-torsion point translations,
identify with abelian scheme isomorphisms $A'\to A''$ up
to involution, and all this in a functorial sense.

\begin{proof} 
Clearly, it suffices to prove that there is an identification 
of groups
$$\isom_{\K_g^{\abs,\pre}}(K',K'')=\Ho_S(S,A''[2])\times\isom_{(\Ag)^{\pre}}(A',A'')$$
which is functorial in $S$.

We first define a group homomorphism
$$\Ho_S(S,A''[2])\times\isom_{\A_g}(A',A'')\xrightarrow{\gamma} \isom_{\K_g^{\abs,\pre}}(K',K'').$$

For any scheme morphism $\alpha$ in $\isom_{Sch}(A',A'')$ and $a$ in  $\Ho_S(S,A''[2])$ we denote by
$\alpha+a$ the isomorphism defined by
$$(\alpha+a)(z):=\alpha(z)+a(v(z))\ \ \ \forall z\in A'.$$

Given any morphism $\varphi$ in $\isom_{\A_g}(A',A'')$
and $a$ in $\Ho_S(S,A''[2])$, the morphism $\varphi+a$ passes to the Kummer quotient
and we denote by $\overline{\varphi + a}$ its class in $\isom_{\K_g^{\abs,\pre}}(K',K'')$. 
The map $\gamma$ sends the pair $(a,\varphi)$
to $\overline{\varphi+a}$.

Let us show that $\gamma$ is surjective.
Consider an isomorphism $\overline\varphi\in\isom_{\K_g^{\abs,\pre}}(K',K'')$ and the associated diagram
$$\xymatrix{
A'\ar@{.>}[r]^{\varphi}\ar[d]_{\pi'}& A''\ar[d]^{\pi''}\\
K'\ar[d]_{u'}\ar[r]^{\overline\varphi}& K''\ar[d]^{u''}\\
S\ar[r]^{\id} & S}$$
where $\varphi$ is an element of $\isom_{Sch}(A',A'')$,
existing by definition of Kummer morphism.
The projection of $2$-torsion points is fixed
by~$\overline\varphi$ because it coincides with the
singular locus of the Kummer schemes. Therefore we have
$$\overline\varphi(\pi'(A'[2]))\subset \pi''(A''[2]).$$
This implies that the codomain of $\varphi\circ \varepsilon$ is
included in~$A''[2]$, or equivalently 
$$\varphi\circ\varepsilon\in\Ho_S(S,A''[2]).$$

Therefore we observe that morphism $\varphi-(\varphi\circ\varepsilon)$
is a lifting of $\overline\varphi$ which preserves the $0$-section $\varepsilon$,
and therefore it preserves the abelian
structure of $A',A''$ by the Rigidity Lemma \cite[p.43]{Mum70}, i.e.~we have $\varphi-(\varphi\circ\varepsilon)\in\isom_{\A_g}(A',A'')$.
This proves that $\gamma$ is surjective.\newline

Observe that the group $\isom_{\A_g}(A'',A'')$ acts on $\isom_{\A_g}(A',A'')$
by composition.
Consider now any element $\overline\varphi'\in\isom_{\K_g^{\abs,\pre}}(K',K'')$,
and a preimage $(a,\varphi')$ of $\overline\varphi'$ via $\gamma$.
Therefore any other preimage of $\overline\varphi'$ is of the form
$(a,\varphi''\circ\varphi)$ with $\varphi''\in\isom_{\A_g}(A'',A'')$. Moreover,
by definition of associated Kummer scheme, $\varphi''$ must be in the subgroup
$\zzz(S)\subset\isom_{\A_g}(A'',A'')$ generated by the natural involution.
Therefore, in consequence also of the equality (\ref{eq_agpre}),
$\gamma$ induces the isomorphism
$$\widetilde\gamma\colon \Ho_S(S,A''[2])\times\isom_{(\Ag)^{\pre}}(A',A'')\to \isom_{\K_g^{\abs,\pre}}(K',K'').$$
\end{proof}

\begin{corollary}
The category fibered in groupoids over schemes $\K_g^{\abs,\pre}$ is a prestack.
\end{corollary}
\begin{proof} Being a prestack is equivalent to the fact that the functor
$\Isom_{\K_g^{\abs,\pre}}(K',K'')$ is representable for any pair
of objects $K'\to S,\ K''\to S$ over the same scheme $S$.
This is a direct consequence of the lemma above, because $\Hoo_S(-,A[2])$  is
clearly represented by the scheme $A[2]$ and given
any two abelian schemes $A',A''$ over $S$, then $\Isom_{\Ag^{\pre}}(A',A'')=\Isom_{\Ag}(A',A'')$
is representable because $\Ag$ is a Deligne-Mumford stack by Theorem \ref{teo_this}.\end{proof}

Consider now a Kummer scheme $K\to S$. There exists a morphism $S\to \K_g^{\abs,\pre}$
canonically associated to $K\to S$. We denote by $A\to S$ an abelian scheme 
whose associated Kummer scheme is $K\to S$.
\begin{lemma}\label{lem_gpre}
With the above notation, we have
$$(\Ag)^{\pre}\times_{\K_g^{\abs,\pre}} S=A[2].$$
In particular the morphism
$G^{\pre}:(\Ag)^{\pre}\to \K_g^{\abs, \pre}$ is representable, finite, étale and surjective.  
\end{lemma}
Therefore if we pullback the prestack morphism $G^{\pre}$ (see diagram (\ref{diagG}))
along $S\to \K_g^{\abs,\pre}$, we obtain a scheme.
Furthermore, even if we used this notation, the pullback is independent
of the choice of the abelian scheme $A$ over $K$, therefore
this lemma implies that the $2$-torsion point scheme $A[2]$ is independent
of this choice too. In the following we will 
denote this finite scheme over $S$ by $K[2]$.
\begin{proof} For any scheme $T$, by the definition of fibered product we have

$$\left((\Ag)^{\pre}\times_{\K_g^{\abs,\pre}}S\right)(T)=\Ho(T,(\Ag)^{\pre}\times_{\K_g^{\abs,\pre}}S)=$$
$$=\left\{\left(\begin{array}{c}
A'\\
\downarrow\\
T\\
\end{array},\ 
T\xrightarrow{v} S,\ \beta\right)\left|\begin{array}{l}
A'\to T\ \m{is an abelian scheme},\\
\beta\in\isom_{\K_g^{\abs,\pre}}(K',v^*K),\\
K'=A'\slash i\m{ is the Kummer scheme associated to } A'.\\ 
\end{array}
\right.\right\}.$$
By Lemma \ref{lemma_aut} we have
$$\isom_{\K_g^{\abs,\pre}}(K',v^*K)=\Ho_T(T,v^*A[2])\times \isom_{(\Ag)^{\pre}}(A',v^*A).$$
Therefore, we have
$$\left((\Ag)^{\pre}\times_{\K_g^{\abs,\pre}}S\right)(T)=$$
$$=\left\{\left(
\begin{array}{c}
A'\\
\downarrow\\
T\end{array},\ T\xrightarrow{v} S,\ \beta_1,\beta_2\right)\left|\begin{array}{l}
A'\to T\ \m{is an abelian scheme},\\
\beta_1\in\Ho_S(T,A[2])\\
\beta_2\in\isom_{(\Ag)^{\pre}}(A',v^*A).\\
\end{array}
\right.\right\}.$$
This concludes the first part of the Lemma. The second assertion follows from the fact that $A[2]\to S$ is finite, étale (because ${\rm char}(k)\neq 2$) and surjective.\end{proof}

We want to prove that $\Ka_g$ is a Deligne-Mumford stack, which unformally speaking
means that there exists an atlas $U\to \Ka_g$ which is a scheme.
In order to prove that, we start by
considering the Deligne-Mumford stack $\A_g$ of principally polarized abelian varieties 
of dimension $g$. By Lemma \ref{lemma_aut} we know that given an abelian variety $A$
and the associated Kummer variety~$K$, then the automorphism group of $K$
is obtained by an extension of $\Aut(A)\slash i$ with a finite set: 
this is  the crucial ingredient in proving that $\K_g^{\abs}$ is a Deligne-Mumford stack.

\begin{theorem}\label{teo2_dm}
The stack $\Ka_g$ is a Deligne-Mumford stack.
Moreover, the natural morphism $\A_g\to \K_g^{\abs}$ factors through a morphism 
$G:\Ag\to \K_g^{\abs}$ which is representable, finite, étale and surjective.  
\end{theorem}
\begin{proof} We have to prove that the diagonal morphism $ \Ka_g\to\Ka_g\times\Ka_g$
is representable, separated and quasi-compact, and moreover that there exists a (representable) étale surjective 
morphism $U\to \Ka_g$.

For the first point, it is sufficient to prove that the diagonal morphism of prestacks
$\K_g^{\abs,\pre}\to\K_g^{\abs,\pre}\times\K_g^{\abs,\pre}$
is representable, separated and quasi-compact. To show that, we consider a morphism $f\colon S\to \K_g^{\abs,\pre}\times\K_g^{\abs,\pre}$
from a scheme $S$. This morphism corresponds to a pair of objects
$K'\to S$ and $K''\to S$ in $\K_g^{\abs,\pre}(S)$, i.e.~two Kummer schemes
over $S$. We fix two abelian schemes $A'$ and $A''$
whose associated Kummer schemes are $K'$ and $K''$ respectively.
For any scheme $T$ we have
$$(S\times_{(\K_g^{\abs,\pre}\times\K_g^{\abs,\pre})}\K_g^{\abs,\pre}) (T)
\Ho(T,S\times_{(\K_g^{\abs,\pre}\times\K_g^{\abs,\pre})}\K_g^{\abs,\pre})$$
$$=\left\{(v\colon T\to S,\beta)|\ \ \beta\in\isom_{\K_g^{\abs,\pre}}(v^*K',v^*K'')\right\}$$
$$=\left\{(v\colon T\to S,\nu)|\ \ \nu\in\Ho_S(T,\Isom_{\K_g^{\abs,\pre}}(K',K'')\right\}$$
$$=\Ho(T,\Isom_{\K_g^{\abs,\pre}}(K',K'')).$$

Therefore $S\times_{(\K_g^{\abs,\pre}\times\K_g^{\abs,\pre})}\K_g^{\abs,\pre}$ 
is represented by the functor $\Isom_{\K_g^{\abs,\pre}}(K',K'')$.
 Using Lemma \ref{lemma_aut} we obtain 
$$\Isom_{\K_g^{\abs,\pre}}(K',K'')=\Hoo_S(-,A''[2])\times\Isom_{\Ag}(A',A'').$$
We observe that, as $\Ag$ is an algebraic stack, $\Isom_{\Ag}(A',A'')$
is represented by a separated and quasi-compact scheme. Moreover, we observe that $\Hoo_S(-,A''[2])$ is represented
by the $S$-scheme $A''[2]$ which is proper over $S$.
As a consequence
$\Isom_{\K_g^{\abs,\pre}}(K',K'')$ is represented by a separated and quasi-compact scheme too.
This proves the first point.\newline

The fact that $\Ag$ is a Deligne-Mumford stack implies that there exists a scheme $U$ and
a representable, étale and surjective morphism $U\to \Ag$.

By Lemma \ref{lem_gpre}, the morphism
of prestacks $G^{\pre}\colon (\Ag)^{\pre}\to\K_g^{\abs,\pre}$ is representable, étale and surjective.
This implies that also the morphism $G:\Ag\to \K_g^{\abs} $ is representable, étale, finite and 
surjective (which also proves the second assertion of the Theorem). Therefore, the composition 
$$U\to \Ag\to \K_g^{\abs}$$
is representable, étale and surjective, and we are done.\end{proof}

\section{Embedded Kummers}\label{emb}

Recall that every principally polarized abelian scheme $\pi\colon A\to S$ of relative dimension $g$ comes with
a globally defined divisor which is fiberwise equivalent to twice
the theta divisor associated to the polarization (see~\cite[Section~VI.2]{GIT}). 
We will denote this divisor by $2\Theta\subset A$.
Since $2\Theta$ is globally generated and has vanishing
higher direct images, we get a morphism
$$\phi_{2\Theta}\colon A \to \p_{\pi}:=\p\left(\pi_*\Oo(2\Theta)\right),$$
where $\p_{\pi}$ is a projective bundle of rank $2^g-1$.
The image of $\phi_{2\Theta}$, which we will denote by $K_{2\Theta}$,
is called the embedded Kummer scheme associated to the abelian scheme $\pi\colon A\to S$.\newline

In this section we work over $\C$, and our goal is to introduce
 a universal Kummer embedding in the case $g=2$,
$$\HH_2\times \C^2\to \p ^{3}\times \p^{3}.$$
From this we will be able to define the notion of embedded Kummer surfaces.

\subsection{Theta functions with characteristics}
We will briefly introduce theta functions with characteristics, which represent a remarkable example of Siegel modular forms and will be very useful in our study of the geometry of embedded Kummer surfaces. 
For a detailed study of theta functions as modular forms see~\cite{ig1}.

Given a point $\tau\in\HH_g$, the factor of automorphy 
\[e(\tau n+q,z)=e^{-\pi i({}^tn\tau n+2{}^tnz)},\ \ n,\,q\in\ZZ^g,\,z\in\CC\]
defines a principal polarization $L_\tau$ on $X_\tau$ (see~\cite[Appendix B]{BL04} for the relationship between line bundles on complex tori and factors of automorphy). 
The unique section up to scalar for $L_\tau$ is defined by the series 
\[\vartheta(\tau,z)=\sum_{n\in\ZZ^g}\mathbb{e}\left(\frac{1}{2}{}^tn\tau n+{}^tnz\right),\]
where $\mathbb{e}(t)=e^{2\pi i t}$. This is called the classical Riemann theta function and it is a holomorphic function in both arguments.
This function cuts out the theta divisor $\T$ 
associated to $(A_\tau,L_\tau)$
\[\T=\{z\in A_{\tau}\mid\vartheta(\tau,z)=0\}.\]
The geometry and the singularities of the theta divisor are extremely interesting in the theory of abelian varieties and have been widely studied (see~\cite{G} for an introduction to the subject). 

For any $x\in A_\tau$, let $t_x\colon A\to A$ be the translation by $x$. Since the first Chern class is invariant by translations, the line bundle $t_x^\ast L_\tau$ defines again a principal polarization on $A_\tau$. 
The sections of the line bundles $t_x^\ast L_\tau$ for $x$ a 2-torsion point are extremely interesting in the theory of modular forms. 
If $m=\ch{a\\b},\;a,\,b\in\{0,1\}^g$, the point $x_m=a\frac{\tau}{2}+b\frac{1}{2}$ is a 2-torsion point in $A_\tau$. The line bundle $t_{x_m}^*L_{\tau}$ 
is symmetric and its unique non-zero section up to scalar is the theta function with characteristic 
\[\vartheta_m(\tau,z)=\tch{a\\b}(\tau,z)=\sum_{n\in\ZZ^g}\mathbb{e}\left(\frac{1}{2}{}^t(n+a/2)\tau (n+a/2)+2{}^t(n+a/2)(z+b/2)\right).\]

The function $\vartheta_m$ is an even function of $z$ if ${}^ta\,b\equiv 0\pmod{2}$ and it is an odd function of $z$ if ${}^ta\,b\equiv 1\pmod{2}$. Correspondingly, the characteristic $m$ is called even or odd. 

Consider the locus
\[\theta_{null}=\left\{\prod_{m\ \text{even}}\vartheta_m(\tau,0)=0\right\}\subset \aa_g,\]
which is an irreducible divisor. It is the locus of those principally polarized abelian varieties for which the theta divisor has a singularity at a 2-torsion point.  
Moreover, for $g=1$, $\theta_{null}$ is empty; for $g=2$ it is equal to the set of reducible abelian varieties;
for $g=3$ it contains the Jacobians of hyperelliptic curves of genus three.

\subsection{Kummer polarization and theta functions: the universal embedded Kummer surface}

Let $(A_\tau,L_\tau)$ be the principally polarized abelian surfact corresponding to a point $\tau$ of~$\HH_2$. The Kummer polarization on $A_\tau$ is then given by the line bundle $L_\tau^2$. A basis of sections of $L_\tau^2$ is given by the second order theta functions 
\[\T[\sigma](\tau,z)=\tch{\sigma\\0}(2\tau,2z),\;\sigma\in\{0,1\}^2.\]

\begin{remark}
We consider the theta group associated to a line bundle
\[\mathcal{G}(L)=\{(x,\varphi)\mid x\in H(L),\,\varphi:L\xrightarrow{\,\simeq\,}t_x^\ast L\}\]
with group law
$(y,\psi)\cdot(x,\varphi)=(x+y,t_x^\ast\psi\circ\varphi)$
for any $(x,\varphi),(y,\psi)\in\mathcal{G}(L)$. 

This group can be entirely described in terms of the type of
the polarization as
\[\mathcal{G}:=\mathcal{G}(2\mathbf{1}_2)=\CC^\ast\times (\ZZ/2\ZZ)^2\times (\ZZ/2\ZZ)^2.\]

Indeed, there exists a natural theta structure, that is a group isomorphism
which is the identity over center $\C^*$,

\begin{equation}\label{natural_theta_structure}
\al_\tau:\mathcal{G}(L_\tau^2)\to\mathcal{G}.
\end{equation}

For the theory of theta structures and other general results about the representation
theory of $\mathcal{G}$ and $\mathcal{G}(L)$, we refer to \cite{Meq1}.
\end{remark}

From these constructions it follows that the induced morphism in projective space is given as
\begin{equation}\label{natural_theta_structure_map}
\begin{split}
\mathcal{F}_{\tau}:A_\tau&\to\PP^3\\
z&\mapsto [\dots,\T[\sigma](\tau,z),\dots].
\end{split}
\end{equation}
The map $\mathcal{F}_{\tau}$ is equivariant for the action of $\mathcal{G}$ on theta functions. 
By~\cite{vG98}, this action is explicitly given by
\begin{equation}\label{action}
(t,x,y)\cdot\T[\sigma]=t(-1)^{(\sigma+x)\cdot y}\T[\sigma+x].
\end{equation}
It then gives an action of the Heisenberg group $\mathcal{H}=\mathcal{G}/\CC^\ast$ on $\PP^3$.
We observe moreover that the Heisenberg group is the maximal group acting on $\PP^3$
and preserving the embedding. To prove this we define an $\mathcal{F}_{\tau}$ automorphism
as a pair $(h',h)$ where $h'\colon A_\tau\to A_\tau$ is a scheme automorphism,
$h\colon \PP^3\to \PP^3$ is a projective automorphism, and $\mathcal{F}_{\tau}\circ h'=h\circ \mathcal{F}_{\tau}$.

\begin{lemma}\label{lem_heis1}
If $(h',h)$ is a $\mathcal{F}_{\tau}$ automorphism, then $h$ is an
element of the Heisenberg group~$\mathcal{H}$.
\end{lemma}
\begin{proof} The fact that any morphism in $\mathcal{H}$
induces an $\mathcal{F}_{\tau}$ automorphism
is implied by the construction above.

For any 
$2$-torsion point $x$ of $A_\tau$, we denote by $(t_x,h_x)$ the $\mathcal{F}_{\tau}$ automorphism
associated to the $x$-translation, and we observe
that $h_x$ is an element of the Heisenberg group.

For any $\mathcal{F}_{\tau}$ automorphism $(h',h)$,
$h$ is in $\mathcal{H}$.
Indeed, consider any $\mathcal{F}_{\tau}$ automorphism $(h',h)$,
by construction $h'$ sends any $2$-torsion point to a $2$-torsion point,
so $h$ fixes $\mathcal{F}_{\tau}(A_\tau[2])$.
If $\varepsilon$ is the
origin of $A_\tau$, then $h'(\varepsilon)=x'$ is a $2$-torsion point of $A_\tau$,
and therefore, taking $(h'\circ t_{-x'},h\circ h_{-x'})$ instead, we can suppose that $h'$ fixes
the origin, $h'(\varepsilon)=\varepsilon$.
Furthermore,
if $x$ is any $2$-torsion point of $A_\tau$, the translation $t_x$ on $A_\tau$ induces
a bijection $t_x^*\colon H^0(L_\tau^2)\to H^0(L_\tau^2)$. The fact that
$h$ fixes $\mathcal{F}_{\tau}(\varepsilon)$ implies that the set 
$$\left\{s\in H^0(L_\tau^2)|\ s(\varepsilon)=0\right\}$$
is fixed, and therefore by taking $t_x^*$ we conclude
that the set of sections with a zero at $x$ is fixed for every $2$-torsion point.
Therefore $h'(x)=x$ for every $2$-torsion
point $x$ of $A_\tau$. Finally, if we fix $\tau$, the second order theta functions $\Theta[\sigma](\tau,z)$ generate
$H^0(L_\tau^2)$. Moreover, for any $\sigma\in \{0,1\}^2$, the zero set of $\Theta[\sigma](\tau,z)$ is included in the set 
of $2$-torsion points.
This means that $h$ acts as the identity on $\PP H^0(L_\tau^2)\cong \PP^3$.
Thus any $h$ as above is of the form $h_x$, i.e.~it is an element of $\mathcal H$.\end{proof}

As a consequence of what we said, the following map is holomorphic,
\begin{align*}
\T_\tau:\HH_2\times\CC^2&\to\PP^3\\
(\tau,z)&\mapsto [\dots,\T[\sigma](\tau,z),\dots].
\end{align*}

We introduce the universal Kummer surface, an object which has been largely studied, see for
example~\cite{vG} and~\cite{RSSS13}.
Consider the map
\begin{equation}\label{univ_Kumm_g}
\begin{split}
\phi:\HH_2\times\CC^2&\to\PP^3\times\PP^3\\
(\tau,z)&\mapsto [\T_\tau(0),\T_\tau(z)].
\end{split}
\end{equation}
The image of this map is a quasi-projective variety of dimension $4$. 
We denote by $\mc W_2$ the closure of the image of $\phi$ inside $\p^3\times \p^3$. Moreover,
we consider the restriction of this morphism to the first components of the products, thus we have
\begin{align*}
\phi_1\colon \H_2&\to\p^3\\
\tau&\mapsto[\T_\tau(0)].
\end{align*}
We denote by $U_2\subset \p V$ the image of $\phi_1$, which is a quasi-projective variety 
(see \cite{vG}). With this notation, we have that 
\begin{equation}\label{family_embedded}
\xymatrix{
\mc W_2\ar@{^{(}->}[r] \ar[d] & U_2\times \p^3\ar[dl]\\
U_2\\}
\end{equation}
\begin{definition}\label{rmk_k}
We call $\mc W_2\to U_2$ the universal embedded Kummer.
We will call embedded Kummer surface any geometric fiber
of the family (\ref{family_embedded}). 
\end{definition}

To every 
abelian surface is naturally associated one embedded Kummer surface, which
is an embedding of the abelian surface in $\PP^3$.
For an indecomposable principally polarized abelian surface $(A,\Theta)$, the associated embedded Kummer surface is
isomorphic to the image of the
canonical embedding $K_{2\Theta}\subset \p H^0(\Oo(2\Theta))$
and to the corresponding abstract Kummer.

In the decomposable case $(A,\Theta)\cong (A_1,\Theta_1)\times(A_2,\Theta_2)$, the image $K_{2\Theta}$
is the reduced scheme associated to the  embedded Kummer surface.\newline

We observe that $U_2$ has exactly one geometric point corresponding to the embedding
of any decomposable abelian surface and such that the geometric fiber over this point is
a double quadric embedded in $\p^3$. Indeed, if $(A,\Theta)$ is a decomposable abelian variety
$A\cong E_1\times E_2$, with a principal polarization $\Theta$, the embedding 
$K_{2\Theta}\subset \p H^0(\Oo(2\Theta))$
is isomorphic to a reduced quadric in $\p^3$.

\section{Moduli space of embedded Kummer surfaces}\label{embsurf}

In this section, we will focus on embedded Kummer surfaces over $\C$
and give a description of the moduli stack and the coarse moduli space of these varieties.

\subsection{Moduli of embedded Kummers}\label{sec_emb}

We start by defining the embedded Kummer functor~$\Ke_2$, and for this we need
to consider the universal embedded Kummer of equation~(\ref{family_embedded}).

\begin{definition}\label{def_ke}
The $\K_2^{\emb, \pre}$ is a fibered category
over schemes whose objects are the diagrams 
$$
\xymatrix{
K\ar@{^{(}->}[r] \ar[d]_{\pi} & \mc B\ar[ld]\\
S,\\}
$$
where $K\to S$ is a flat morphism of schemes, $\mc B$
is a projective bundle of rank $3$ over~$S$,
and every geometric fiber of this diagram
is an embedded Kummer surface (see Definition \ref{rmk_k}).

Given a pair of embedded Kummer families
$S\leftarrow K\hookrightarrow \mc B$ and $S'\leftarrow K'\hookrightarrow \mc B'$,
a morphism between them in $\K_2^{\emb,\pre}$
is a pair of scheme morphisms 
$(u\colon S\to S',\ v\colon \mc B\to \mc B')$,
 such that
in the following diagram every square is cartesian
$$\xymatrix{
K\ar[d]_{\pi}\ar[dr]\ar@{^{(}->}[r] & \mc B\ar[rd]^{v}\\
S\ar[dr]_{u} & K'\ar@{^{(}->}[r] \ar[d]_{\pi'}&\mc B'\\
& S'.\\}
$$

By construction $\K_2^{\emb,\pre}$ is a prestack. The stack $\Ke_2$ 
of embedded Kummer varieties is the stack associated to $\K_2^{\emb,\pre}$.
\end{definition}

\begin{proposition}\label{prop_g2}
The universal embedded Kummer is a flat family.
\end{proposition}
\begin{proof} As shown in \cite{vG}, the closure $\overline{\mc W}_2$
is cut out by a single equation and this implies the flatness over $U_2$.
\end{proof}

We denote by
$$p(t):=2t^2+2,$$
the Hilbert polynomial of any degree $4$ hypersurface inside $\p^3$,
which therefore is the Hilbert polynomial of the general embedded Kummer surface. 
By the universal property of the Hilbert scheme,
we have a morphism
$$\phi'\colon U_2\to \hilb_{\p^3}^{p(t)}.$$

\begin{remark}
The lack of the flatness property is the obstruction
to generalize our results to $g>2$. Indeed, the universal embedded Kummer
is not proven to be flat for greater $g$,
and therefore we don't have an immersion of the general $U_g$ scheme
in the Hilbert scheme.
\end{remark}


We recall from~\cite{vG98} that the group $\Gamma(2,4)$ is the subgroup of the integral symplectic group of rank $4$
such that its action on $\HH_2$ preserves the isomorphism class of a principally polarized abelian
variety plus a theta structure. The quotient $\HH_2\slash\Gamma(2,4)$ is denoted
$\aa_2(2,4)$ and is the coarse moduli space of principally polarized abelian varieties with a theta structure.

\begin{lemma}\label{lemma_factor}
The morphism $\phi'\circ\phi_1\colon \HH_2\to\hilb_{\p^3}^{p(t)}$
factors through $\aa_2(2,4)$.
\end{lemma}
\begin{proof} It suffices to prove that two points of $\HH_2$ corresponding
to the same theta structure on a certain abelian surface $(A,2\Theta)$, are
sent in the same point of $\hilb_{\p^3}^{p(t)}$. In fact, given
a theta structure on the abelian surface above, the morphism
$A\to \PP^3$ is uniquely determined.
This implies that the two images inside $\p^3$ coincide, and therefore
they correspond to the same point of $\hilb_{\p^3}^{p(t)}$.\end{proof}

As another consequence of the above setting, 
the action of $\hh$ on $\p^3$ induces an action
of the same group on $\hilb_{\p^3}^{p(t)}$.
If we denote by $\widetilde U_2$ the image of $\phi'$ inside $\hilb_{\p^3}^{p(t)}$, then
$\widetilde U_2$ is fixed by the $\hh$ action.
Indeed, the action of $\hh$ on $\hilb_{\p^3}^{p(t)}$
changes the coordinates of $\p^3$ and thus sends any embedded Kummer  variety
to another embedded Kummer  variety. 

\begin{theorem}
The moduli stack $\Ke_2$ is the quotient stack
$[\widetilde U_2\slash\hh]$, and therefore it is a Deligne-Mumford stack.
\end{theorem}
\begin{proof} To prove this, we first define a morphism of categories fibered in grupoids
$$\Psi^{\pre}\colon\K^{\emb,\pre}_2\to[\widetilde U_2\slash \hh].$$
As $[\widetilde U_2\slash \hh]$ is a stack, this induces a natural morphism $\Psi\colon \Ke_2\to[\widetilde U_2\slash \hh]$.
Finally, we will prove that $\Psi$ is an isomorphism.

Consider an object of $\K^{\emb,\pre}_2$, i.e.~a family of embedded Kummer  varieties
$$\xymatrix{ K\ar[d]_\pi \ar@{^{(}->}[r] &\mc B \ar[dl]\\ S.}$$
By definition, every geometric fiber of $\pi$ is an embedded Kummer variety,
therefore by Lemma~\ref{lem_heis1},
$\mc B$ is naturally a
$\hh$-bundle over $S$.
We denote by $\xi\colon E\to S$ the principal $\hh$-bundle associated to $\mc B$. 
We observe that $\xi^*\mc B$ is a trivial projective bundle, i.e.~we have
$$\xymatrix{
E\times_S K\ar[d] \ar@{^{(}->}[r] & E\times \p V\ar[dl]\\
E\\ }$$
for some vector space $V$.
By the universal property of the Hilbert scheme, there exists
a unique morphism $m\colon E\to \hilb_{\p V}^{p(t)}$
which is $\hh$-equivariant. By construction,
 we have that the image of $m$ is included in $\widetilde U_2$,
 and finally the object $\Psi^{\pre}(K\to S)$ is
$$\xymatrix{
E\ar[r]^{m} \ar[d]_{\xi} & \widetilde U_2\\
S.
}
$$

We construct now the image via $\Psi$ of a morphism $(u, v)$ in $\K^{\emb,\pre}_2$ between two families of embedded Kummer varieties
$S'\leftarrow K'\hookrightarrow \mc B'$ and $S''\leftarrow K''\hookrightarrow \mc B''$. 
We denote by $E'\to S'$ and $E''\to S''$ the principal
$\hh$-bundles associated to $\mc B'$ and $\mc B''$ respectively.
Moreover, we denote by $m'\colon E'\to \widetilde U_2$ and $m''\colon E''\to \widetilde U_2$
the associated $\hh$-equivariant morphisms.
Then
the morphism $v\colon \mc B'\to \mc B''$ induces
and $\hh$-equivariant morphisms $E'\to E''$ commuting
with $m'$ and $m''$, and such that the following diagram is
cartesian
$$
\xymatrix{
E'\ar[r]\ar[d]& E''\ar[d]\\
S'\ar[r] & S''.\\
}
$$
This map between $E'$ and $E''$ is exactly the image morphism $\Psi(u,v)$.

To prove that $\Psi$ is a morphism of categories it suffices to prove that given a morphism
$(u',v')$ from $S'\leftarrow K'\hookrightarrow B'$ to $S''\leftarrow K''\hookrightarrow \mc B''$ and a morphism $(u'',v'')$ from 
$S''\leftarrow K''\hookrightarrow \mc B''$ to $S'''\leftarrow K'''\hookrightarrow \mc B'''$,
then $\Psi\left((u'',v'')\circ(u',v')\right)=\Psi(u'',v'')\circ\Psi(u',v')$. This follows
easily from the construction above.\newline

It remains to prove that $\Psi$ is an isomorphism. It suffices to show that for every
scheme~$S$, the functor $\Psi_S\colon \Ke_2(S)\to[\widetilde U_2\slash\hh](S)$ is fully faithful 
and essentially surjective, i.e.~it is an equivalence of categories (see \cite[Lemma 5.1 \S12]{acgh11}).

Fully faithfulness means that there is a bijection between the set of automorphisms
of a family $S\leftarrow K\hookrightarrow \mc B$ and the set of automorphisms of the associated $\hh$-torsor:
$$\Ho_{\Ke_2}(K,K)\cong\Ho_{[\widetilde U_2\slash \hh]}(\Psi(K),\Psi(K)).$$
This is equivalent to proving the same property for the prestack $\K_2^{\emb,\pre}$ and 
morphism $\Psi$,
and it follows from the definitions.

Essential surjectivity means that any object of $[\widetilde U_2\slash \hh]$
is isomorphic to the image via $\Psi$ of some object in $\Ke_2$. This is equivalent to say
that for any object of $E\to S$ of $[\widetilde U_2\slash \hh]$, there exists
a covering $\{S_i\}$ of the scheme $S$ such that every restricted object
$E|_{S_i}\to S_i$ is isomorphic to the image via $\Psi$ of some object in $\K^{\emb,\pre}_2$.
To prove this,
we consider the universal family
$$\xymatrix{
\mc Y \ar[d] \ar@{^{(}->}[r] &\widetilde U_2\times\p V\ar[ld]\\
\widetilde U_2\\
}$$
which is the restriction of the universal family over the Hilbert scheme.
We consider an object of $[\widetilde U_2\slash \hh]$, i.e.~an $\hh$-torsor plus an $\hh$-equivariant morphism to $\widetilde U_2$
$$\xymatrix{
E \ar[r]^m\ar[d] & \widetilde U_2\\
S.\\
}$$
The group $\hh$ acts equivariantly and freely on $E$ and $\mc Y$, so we consider a covering
$\{S_i\}$ of $S$ that induces a trivialization of $E$, equivalently such that
$E|_{S_i}\cong S_i\times \hh$ for any $S_i$.
We consider now the quotient family
$$
\xymatrix{
(m^*\mc Y|_{S_i})\slash \hh\ar@{^{(}->}[r]\ar[d]_{\pi_i} & S_i\times \p V\ar[ld]\\
S_i.\\
}
$$

By the definition of $\mc Y$, for any $S_i$ this is a family of embedded
varieties whose Hilbert polynomial is $p(t)$, i.e.~a family
of embedded Kummers, therefore an object of $\K_2^{\emb,\pre}(S_i)$.
This completes the proof.\end{proof}

In order to give a description of the coarse moduli space of embedded Kummer surfaces,  
we define a $1$-cycle numerical class on the coarse moduli space
of abelian surfaces $\aa_2$.
\begin{definition}
Given a $1$-dimensional abelian variety $B$, we denote by $C_1^{(B)}\subset\aa_2$
the locus of abelian surfaces isomorphic to $B\times E$, where $E$ is any $1$-dimensional abelian variety.
\end{definition}

\begin{remark}\label{rmk_xi}
The curves $C_1^{(B)}$ in $\aa_2$ are all linearly equivalent, independently of the 
abelian variety $B$. Indeed, there exists a morphism 
$$\Xi\colon\aa_1\times\aa_1\to \aa_2,$$
sending any geometric point $[B]\times [E]$ of the product to the 
point associated to the abelian surface $[B\times E]$. Therefore
$C_1^{(B)}=\Xi\left(\{[B]\}\times \A_1\right)$,
and as $\aa_1\cong \p^1$, the
curve $C_1^{(B)}$ is in the same equivalence class $C_1$ for any $B$.
\end{remark}

We define now a natural map $\beta_2\colon\A_2\to \Ke_2$.
Consider any abelian scheme $A\to S$ in $\A_2$. By definition
there exists an associated morphism $p_A\colon S\to\aa_2$
Consider a lifting $\bar p_A\colon S\to \HH_2$, and the morphism $\phi_1\colon \HH_2\to U_2$.
By pullbacking the universal family $\mc W_2$ over $S$ via $\phi_1\circ \bar p_A$, we obtain
a family of embedded Kummer varieties
$$
\xymatrix{
K\ar@{^{(}->}[r]\ar[d] & \PP V\times S\ar[dl]\\
S.\\
}
$$
Changing the lifting $\bar p_A$ changes the immersion morphism
of $K$ by a projective automorphism of $\PP V$. Therefore
the isomorphism class of the family is well defined,
and this is the image of $A\to S$ in $\Ke_2$ via $\beta_2$.

This induces, passing to the associated coarse moduli spaces, a map $\aa_2\to \kk_2$,
also called~$\beta_2$. We denote by
$\aa_{1,1}$ the image of the map $\Xi$ in Remark \ref{rmk_xi}, i.e.~the sublocus
of decomposable abelian surfaces.

\begin{theorem}
The map $\beta_2\colon\aa_2\to \kk_2$ is the contraction of the class $C_1$,
the equivalence class of $C_1^{(B)}$ for any $B$.
\end{theorem}
\begin{proof} 
For every geometric point $\bar s$ of $\aa_{1,1}$, $\beta_2(\bar s)$
is always the same point $y\in \kk_2$ corresponding to the double quadric embedded in $\p^3$.
We observe that $\aa_2\backslash\aa_{1,1}$ is isomorphic to $\kk_2\backslash\{y\}$,
because families of irreducible abstract Kummer surfaces are in natural bijection
with families of irreducible embedded Kummer surfaces.
Therefore every point in $\aa_{1,1}$ is contracted to $y$ and every point outside $\aa_{1,1}$ is not contracted.
As $\aa_{1,1}$ is covered by curves linearly equivalent to $C_1$, this proves that $\beta$
is exactly the contraction of $C_1$.\end{proof}

\begin{remark}
We observe that since no point outside $\aa_{1,1}$ is contracted, the proof above implies that any curve
not contained in $\aa_{1,1}$ is not equivalent to $C_1$. Equivalently, $\aa_{1,1}$ is the locus covered
by all the curves in the class $C_1$.
\end{remark}

We have thus proven that  the moduli space $\Ke_2$ is a Deligne-Mumford stack
and its coarse moduli space $\kk_2$ is obtained by contracting a particular curve class inside~$\aa_2$.
It remains an open question if such a theorem is true more generally for any $g$. We observe
that if the coarse moduli space $\kk_g$ exists, abelian varieties
of the type $B_{g-1}\times E$, where $B_{g-1}$ is a fixed ($g-1$)-dimensional abelian variety and
$E$ any elliptic curve, run a curve $C_1^{(B)}$ always in the same class $C_1$,
and they are represented by the same geometric point of~$\kk_g$.
We conjecture that $\kk_g$ is obtained via the contraction of the curve class $C_1$ on~$\aa_g$
for any $g>2$. The key obstruction to this remains the proof that $\Ke_g$ is Deligne-Mumford.

In \cite{SB16}, Shepherd-Barron treat the locus $\aa_{1,g-1}^P\subset \aa_g^P$ of decomposable 
abelian varieties of the type $B_{g-1}\times E$
inside the perfect compactification of $\aa_g$. He proves that for $g\leq 11$,
  $\aa_{1,g-1}^P$ is the exceptional locus of a semi-ample line bundle and therefore there exists a contraction
$\aa_g^P\to V_g$ of the class $C_1$. This algebraic variety $V_g$
would contain $\kk_g$ as a dense a open subset if our
conjecture is correct.

Furthermore, in \cite[Corollary 3.4]{S-B} the same author shows that in the case $g=11$,
our conjecture would imply
that a suitable compactification of $\kk_{11}$ is the canonical model of~$\aa_{11}$.

\subsection{Equations for a family of embedded Kummer surfaces}

The universal Kummer variety in dimension 2 can be described in equations, as
done for example in \cite[Example 1.1]{RSSS13}. In this section we show these equations
and show that they also parametrize the embedded Kummer surfaces of reducible abelian varieties and count them with the right multiplicity.

If $g=2$ the map~\eqref{univ_Kumm_g} is the map 
\begin{equation*}
\begin{split}
\phi:\HH_2\times\CC^2&\to\PP^{3}\times\PP^{3}\\
(\tau,z)&\mapsto [\T_\tau(0),\T_\tau(z)].
\end{split}
\end{equation*}
Let $x_\sigma=\Theta[\sigma](\tau,z)$ and $u_\sigma=\Theta[\sigma](\tau,0)$. Then set $\mathbf{u}=(u_{00},u_{01},u_{10},u_{11})$ and $\mathbf{x}=(x_{00},x_{01},x_{10},x_{11})$ and take $(\mathbf{u},\mathbf{x})$ as the coordinates in the above $\PP^{3}\times\PP^{3}$. 

By~\cite{vG}, we know that the closure of the image of $\phi$ is given by a single equation $F(\mathbf{u},\mathbf{x})$ which involves Heisenberg invariant polynomials.

By \eqref{action}, it is easy to see that the action of the Heisenberg group on the $\mathbf{x}$ coordinates is generated by the following transformations
\begin{align*}
(x_{00},x_{01},x_{10},x_{11})&\mapsto (x_{10},x_{11},x_{00},x_{01}),\\
(x_{00},x_{01},x_{10},x_{11})&\mapsto (x_{01},x_{00},x_{11},x_{10}),\\
(x_{00},x_{01},x_{10},x_{11})&\mapsto (x_{00},x_{01},-x_{10},-x_{11}),\\
(x_{00},x_{01},x_{10},x_{11})&\mapsto (x_{00},-x_{01},x_{10},-x_{11}).
\end{align*}
The space of Heisenberg-invariant quartics is generated by the following ones:
\begin{align*}
P_0&=x_{00}^4+x_{01}^4+x_{10}^4+x_{11}^4;\\
P_1&=2(x_{00}^2x_{01}^2+x_{10}^2x_{11}^2);\\
P_2&=2(x_{00}^2x_{10}^2+x_{01}^2x_{11}^2);\\
P_3&=2(x_{00}^2x_{11}^2+x_{01}^2x_{10}^2);\\
P_4&=4x_{00}x_{01}x_{10}x_{11}.
\end{align*}
The equation of the universal Kummer variety is then given as:
\begin{equation}\label{F}
F(\mathbf{u},\mathbf{x})=\det
\begin{pmatrix}
P_0&P_1&P_2&P_3&P_4\\
u_{00}^3&u_{00}u_{01}^2&u_{00}u_{10}^2&u_{00}u_{11}^2&u_{01}u_{10}u_{11}\\
u_{01}^3&u_{00}^2u_{01}&u_{01}u_{11}^2&u_{01}u_{10}^2&u_{00}u_{10}u_{11}\\
u_{10}^3&u_{10}u_{11}^2&u_{00}^2u_{10}&u_{01}^2u_{10}^2&u_{00}u_{01}u_{11}\\
u_{11}^3&u_{10}^2u_{11}&u_{01}^2u_{11}&u_{00}^2u_{11}&u_{00}u_{01}u_{10}
\end{pmatrix}.
\end{equation}

As we already mentioned, the locus of decomposable principally polarized abelian surfaces is the divisor
\[\theta_{null}=\left\{\prod_{m\ \text{even}}\vartheta_m(\tau)=0\right\}\subset \mathcal{A}_2,\]
where $\vartheta_m(\tau):=\vartheta_m(\tau,0)$.
Whereas $\theta_{null}$ is irreducible in $\mathcal{A}_2$, its lifting
to $\H_2$
has 10 irreducible components, each corresponding to the vanishing of one of the even theta constants $\vartheta_m(\tau)$.

Riemann's addition formula gives a way to find equations for these components that are particularly useful to our purpose. Indeed, these equations can be written in terms of second order theta constants which do not satisfy any algebraic relation in degree 2. 
For $\e,\,\de\in\{0,1\}^g$, one has the following formula:
\begin{equation}\label{addition}
\tch{\e\\\de}(\tau)^2=\sum_{\sigma\in\{0,1\}^2}(-1)^{\de\cdot\sigma}\Theta[\sigma](\tau,0)\Theta[\sigma+\e](\tau,0).
\end{equation}
We write a theta characteristic in rows, i.e. $m=\ch{\e_1&\e_2\\\de_1&\de_2}$. We will order the 10 even characteristics in the following way:
\begin{align*}
m_1=\ch{0&0\\0&0},\,
m_2=\ch{0&0\\0&1},\,
m_3=\ch{0&0\\1&0},\,
m_4=\ch{0&0\\1&1},\,
m_5=\ch{0&1\\0&0},\,\\
m_6=\ch{0&1\\1&0},\,
m_7=\ch{1&0\\0&0},\,
m_8=\ch{1&0\\0&1},\,
m_9=\ch{1&1\\0&0},\,
m_{10}=\ch{1&1\\1&1}.
\end{align*}
Correspondingly, we can easily find the equations of the irreducible components of $\theta_{null}$ by~\eqref{addition}. These are the zero locus of the 10 quadrics:
\begin{align*}
q_1(\mathbf{u})=&\; u_{00}^2+u_{01}^2+u_{10}^2+u_{11}^2,\\ 
q_2(\mathbf{u})=&\; u_{00}^2-u_{01}^2+u_{10}^2-u_{11}^2,\\
q_3(\mathbf{u})=&\; u_{00}^2+u_{01}^2-u_{10}^2-u_{11}^2, \\
q_4(\mathbf{u})=&\; u_{00}^2-u_{01}^2-u_{10}^2+u_{11}^2,\\
q_5(\mathbf{u})=&\; 2\,(u_{00} u_{01}+u_{10} u_{11}),\\
q_6(\mathbf{u})=&\; 2\,(u_{00} u_{01}-u_{10} u_{11}),\\
q_7(\mathbf{u})=&\; 2\,(u_{00} u_{10}+u_{01} u_{11}),\\
q_8(\mathbf{u})=&\; 2\,(u_{00} u_{10}-u_{01} u_{11}),\\
q_9(\mathbf{u})=&\; 2\,(u_{01} u_{10}+u_{00} u_{11}),\\
q_{10}(\mathbf{u})=&\; 2\,(u_{00} u_{11}-u_{01} u_{10}).
\end{align*}

By computations in Macaulay2~\cite{M2}, we compute equation~\eqref{F} on the irreducible components of the $\theta_{null}$ divisor. Let $k_i=F(\mathbf{u},\mathbf{x})|_{\{q_i(\mathbf{u})=0\}}$, then  
\[k_i(\mathbf{u},\mathbf{x})=p_i(\mathbf{u})\,f_i(\mathbf{x})^2,\quad i=1,\dots,10,\]
where $p_i(\mathbf{u})$ is a polynomial in $u_{00},\,u_{01},\,u_{10},\,u_{11}$ and 
\begin{align*}
f_1(\mathbf{x})=& x_{00}^2+x_{01}^2+x_{10}^2+x_{11}^2,\\ 
f_2(\mathbf{x})=& x_{00}^2-x_{01}^2+x_{10}^2-x_{11}^2,\\
f_3(\mathbf{x})=& x_{00}^2+x_{01}^2-x_{10}^2-x_{11}^2, \\
f_4(\mathbf{x})=& x_{00}^2-x_{01}^2-x_{10}^2+x_{11}^2,\\
f_5(\mathbf{x})=& x_{00} x_{01}+x_{10} x_{11},\\
f_6(\mathbf{x})=& x_{00} x_{01}-x_{10} x_{11},\\
f_7(\mathbf{x})=& x_{00} x_{10}+x_{01} x_{11},\\
f_8(\mathbf{x})=& x_{00} x_{10}-x_{01} x_{11},\\
f_9(\mathbf{x})=& x_{01} x_{10}+x_{00} x_{11},\\
f_{10}(\mathbf{x})=&x_{00} x_{11}-x_{01} x_{10}.
\end{align*}
It is interesting to see that we get essentially the same quadrics but doubled. This mirrors the fact that, moving from an irreducible abelian surface to a decomposable one, the degree of the image of the map defined by the line bundle $L_\tau^2$ drops. Then, a flat family of embedded Kummer surfaces should carry a suitable non-reduced structure.

\section{Toward compactifications: torus rank 1 degenerations of principally polarized abelian varieties}\label{degen}

Regarding the moduli space of principally polarized abelian varieties over $\C$, of dimension~$g$, there are essentially two kinds of compactifications that have been studied: Satake's compactifications and toroidal compactifications. 

Satake's compactification is developped the setting of bounded symmetric domains:
it is the projective scheme associated to the graded ring of scalar modular forms with respect to the full modular group.
 As a set
\[\overline{\bm A}_g^{Sat}=\aa_g\sqcup\aa_{g-1}\sqcup\dots\sqcup\aa_0,\]
therefore it has a boundary of $\overline{\bm A}_g^{Sat}$ has high codimension. Moreover the Satake's compactification of $\aa_g$ is highly singular along the boundary.
We will denote by $\overline{\bm A}_g^{Sat}$ the Satake's compactification of $\aa_g$.

The toroidal compactifications are a class of different compactifications of $\aa_g$,
each one corresponding to a way of decomposing the cone of real positive quadratic forms in
$g$ variables. They lack some of the bad properties of $\overline{\bm A}_g^{Sat}$,
in particular the boundary of any toroidal compactification is a divisor.

The investigation of degenerations of abelian varieties is a crucial point in the study of toroidal compactifications of $\aa_g$. These indeed represent the boundary points of such compactification. Only recently the knowledge of deeper boundary strata of the toroidal compactification  corresponding to the perfect cone decomposition has been increased (see~\cite{GH2011}). 

The so-called Mumford's partial compactification of $\aa_g$, denoted by $\overline{\bm A}_g^{Mum}$, is contained in any toroidal compactification and we can think of it as representing the most superficial strata of the boundary of any toroidal compactification. 
This compactification can be seen as the blow-up of the open set $\aa_g\sqcup\aa_{g-1}$ in $\aa_g^{Sat}$ along its boundary~$\aa_{g-1}$. 

As explained in~\cite{Mum83}, $\overline{\bm A}_g^{Mum}$ is the coarse moduli space of principally polarized abelian varieties of dimension $g$ and their torus rank 1 degenerations,  which represent the points on the boundary of this partial compactification.

Following~\cite{Mum83} we give the definition of a torus rank 1 degeneration of dimension $g$.
We recall that if $B$ is a principally polarized abelian variety and we denote by $\hat B$ its dual, then $B$ is canonically
identified with $\hat B$ thanks to the principal polarization. Moreover, consider an extension $G$ of $B$ by $\CC^*$,
that is an algebraic group who fits into the short exact sequence
$$0\to \CC^*\to G\to B\to 0.$$
The group of $B$ extensions by $\CC^*$
is isomorphic to $\hat B$,
$$\ext^1(B,\CC^*)\cong \hat B=B.$$

\begin{definition}\label{def_glue}
Take $B$ a principally polarized abelian variety of dimension $g-1$,
let $G$ be an algebraic group which is an extension of $B$ by $\CC^\ast$
and suppose the class of $G$ in $\ext^1(B,\CC^*)$ is identified with the point $b\in B$.
Consider $G$ as a $\CC^\ast$-bundle over $B$ and let $\widetilde{G}$ be the associated $\PP^1$-bundle. 
Then $\widetilde{G}-G$ equals $\widetilde{G}_0\sqcup\widetilde G_\infty$, the union of two sections of $\widetilde G$ over $B$;
 $\overline{G}$ is the variety obtained by gluing $\widetilde{G}_0,\,\widetilde G_\infty$ with a translation by $b$,
then $\overline G$ is a complete $g$-dimensional variety.

Consider an ample Cartier divisor $D$ on $\overline{G}$ such that $h^0(D)=1$, not equivalent
to the divisor of the glued sections $\widetilde G_0,\widetilde G_\infty$. Then
$(\overline{G},D)$ is a torus rank 1 degeneration.
Such a degeneration will be denoted by $X(B,b)$. A variety constructed in such a way is also known as a semi-abelic variety of torus rank 1.
\end{definition}

By~\cite{Mum83}, we know also that the complete variety $\overline{G}$ is the limit of a $g$-dimensional abelian varieties. 
By the equation~\eqref{moduli_ppav} on $\aa_g$, one can give an analytic method to construct this limit. This will also give us a way to write explicitly the divisor $D$ and hence we get an analytic description of a torus rank 1 degeneration of an abelian variety.

Recall that any point $\tau\in\HH_g$ defines the principally polarized abelian variety
\[A_\tau=\CC^g/\tau\ZZ^g+\ZZ^g\]
with theta divisor 
\[\T=\{z\in A_{\tau}\mid\vartheta(\tau,z)=0\},\]
where $\vartheta$ is the theta function with $0$ characteristic. 
Then a point on the boundary of  $\overline{\bm A}_g^{Mum}$ gives a degeneration of the abelian variety $A_{\tau}$ 
when $\op{Im}\tau_{11}\to\infty$ and $\tau_{ij}\to\tau_{ij}'$ finite limite,  for $i\ge2$ or $j\ge2$. If 
\[\tau\to\begin{pmatrix}
i\infty & {}^tb\\
b &\tau'
\end{pmatrix},
\]
we get the semi-abelic variety $X(A_{\tau'},b)$ as the limit. 
By Definition \ref{def_glue}, the $0$-section and the $\infty$-section of the $\PP^1$-bundle over $A_{\tau'}$ are glued compatibly
with the shift by the point $b\in A_{\tau'}$. Formally, the identification is
\[(z',0)\sim(z'-b,\infty).\] 

We can also compute the degeneration of the theta divisor by analytic methods. 
Following~\cite[formula (4)]{GH2011}, we define the semi-abelic theta function to be
\[\lim_{\tau_{11}\to i\infty}\vartheta\left(\sm{\tau_{11} & {}^tb\\
b &\tau'{}},\sm{z_1-\tau_{11}/2\\z'}\right)=\vartheta(\tau',z')+q\,\vartheta(\tau',z'+b),\]
where $z'=(z_2,\dots,z_g)$ is the analytic coordinate on $B$ and $q:=e^{2\pi iz_1}$. 
Note that $q$ can be viewed as the coordinate in the fiber $\CC^\ast$ of $G$ and the theta divisor $\Xi$ on $A_{\tau'}$ as the zero locus of $\vartheta(\tau',z')$. 
The degeneration of the theta divisor is the zero locus of the semi-abelic theta function.

In~\cite{GH2011} one also finds a description of the ``2-torsion points'' on a torus rank 1 degeneration of a principally polarized abelian variety as fixed points of an involution that preserves the semi-abelic theta divisor. 

Let $X=X(B,b)$ be a torus rank 1 degeneration of a principally polarized abelian variety constructed as above. Following~\cite{GH2011} we will first describe a uniformization of the smooth part of $X$. 
Let $B=\CC^{g-1}/(\ZZ^{g-1}\tau\oplus\ZZ^{g-1})$ for $\tau\in\HH_{g-1}$. The semi-abelian variety corresponding to a point $b\in B$ is the quotient of the trivial $\CC^\ast$-bundle $\CC^{g-1}\times\CC^\ast\to\CC^{g-1}$ by the group $\ZZ^{2(g-1)}\cong\ZZ^{g-1}\tau\oplus\ZZ^{g-1}$ acting via
\[(n,m):(z,x)\mapsto(z+n\tau+m,\mathbb{e}(-{}^tn\,b)x),\]
where as before $\mathbb{e}(t)=e^{2\pi i t}$.

Then the involution $j$ on the semi-abelic variety that preserves the semi-abelic theta divisor is given as follows:
\[j(z,x)=(-(z+b),x^{-1}).\]
Note that the involution $j$ exchanges the $0$-section and the $\infty$-section. 

The fixed points of the involution $j$ are the limits of the $2$-torsion points on the degenerating family of abelian varieties. There are two possibilities: they can lie either on the smooth part of $X$ or on its singular locus.
The points
\begin{align*}
z&=-(\tau\e+\de+b), \e,\,\de\in\{0,1\}^{g-1},\\
x&=\pm \mathbb{e}(-{}^t\e\,b).
\end{align*}
are the $2^{2g-1}$ fixed points of $j$ on the smooth locus of~$X$. 
Each point is the limit of a single family of $2$-torsion points on smooth abelian varieties. 

On the other hand, the fixed points of the involution $j$ lying on the singular locus
are of the form
\[z=-(\tau\e+\de), \e,\,\de\in\{0,1\}^{g-1}.\]
Each point must be counted with multiplicity $2$ since it is the limit of two $2$-torsion points on a smooth abelian variety. So we get $2^{g-1}$ fixed points of the involution on the singular locus of~$X$. 

Note that on a semi-abelic variety it does not make sense to say that a fixed point of the involution is even or odd. 
Indeed, a fixed point of $j$, say on the singular locus, can be the limit of odd $2$-torsion points and of even $2$-torsion points (see~\cite{GH2011} for the case of elliptic curves).  

We would like to understand if the 
computation of the degeneration of the map
\begin{align*}
\phi:\HH_g\times\CC^g&\to\PP V\times\PP V\\
(\tau,z)&\mapsto [\T_\tau(0),\T_\tau(z)].
\end{align*}
defining the universal Kummer variety  might lead to the construction of a compactification for the moduli space of embedded Kummer varieties. We leave it here as an open question.

One could try to compute the degeneration of the coordinates $\T_\tau(z)$ and $\T_\tau(0)$ separately. 
Regarding $\T_\tau(z)$, we have seen that they define the embedding of the embedded Kummer varieties of a principally polarized abelian variety. 
For $A_\tau$ a principally polarized abelian variety, let 
\begin{align*}
f:A_\tau&\to\PP^{2^g-1}\\
\tau&\mapsto[\T_\tau(z)]
\end{align*}
be the map defined by second order theta functions. By~\cite{Mum83} we know that such a map extends to torus rank 1 degenerations. If $\op{Im}\tau_{11}\to\infty$, it is defined by the $2^g$ ``theta functions"
\begin{align*}
&\T[\e](\tau',z')+q^2\T[\e](\tau',z'+b),\\
&q\T[\e](\tau',z'+b/2),
\end{align*} 
for $\e\in\{0,1\}^{g-1}$.
In order to compute such degenerations, one can perform the Fourier-Jacobi expansion. By direct computations one has that
\begin{itemize}
\item $\T[0\,\e](\tau,z)$ degenerates to $\T[\e](\tau',z')+q^2\T[\e](\tau',z'+b)$,
\item $\T[1\,\e](\tau,z)$ degenerates to $q\T[\e](\tau',z'+b/2)$.
\end{itemize}
So, at least on Mumford's partial compactification, this this would give the degeneration of the coordinates $\T_\tau(z)$.

What about the coordinates $\T_\tau(0)$? 
If we denote by $\aa_g(2,4)$ the moduli space of principally polarized abelian varieties with theta structure, we know that they define a map
\begin{align*}
\Phi:\aa_g(2,4)&\to\PP^{2^g-1}\\
\tau&\mapsto[\T_\tau(0)]
\end{align*}
such that the closure of the image of $\Phi$ is the Satake compactification of $\aa_g(2,4)$.

It would be interesting to see if there exists a toroidal compactification $\aa_g(2,4)^{tor}$ of $\aa_g$ for which it is possible to define a projective map
\[\Psi:\aa_g(2,4)^{tor}\to\PP^n\]
for some $n$ such that the restriction of $\Psi$ to $\aa_g(2,4)$ coincides with $\Phi$. We expect that this map should be given by some ``degeneration'' of theta constants.

\section*{Acknowledgments}
We are grateful to Alfio Ragusa, Francesco Russo and Giuseppe Zappal\`a
for organizing Pragmatic 2015, where this work started. Special thanks are also due to Giulio Codogni and Filippo Viviani for introducing us to the problem and for their support and comments.
Finally, many thanks to the anonymous reviewer of the
first version of this paper, for his precious suggestions.


\begin{thebibliography}{0} 
\bibitem{acgh11} E. Arbarello, M. Cornalba, P. A. Griffiths: Geometry of Algebraic Curves Volume II, 
 Grundlehren der Mathematischen Wissenschaften, 268. Springer, Heidelberg, 2011. With a contribution by J. D. Harris.
 
 \bibitem{ACV03}
Dan Abramovich, Alessio Corti, and Angelo Vistoli.
\newblock Twisted bundles and admissible covers.
\newblock {\em Comm. Algebra}, 31(8):3547--3618, 2003.
\newblock Special issue in honor of Steven L. Kleiman.
 
\bibitem{BL04} C. Birkenhake, H. Lange: Complex abelian varieties. Second edition. Grundlehren der Mathematischen Wissenschaften, 302. Springer-Verlag, Berlin, 2004.

\bibitem{FC} G. Faltings, C-L. Chai: Degeneration of abelian varieties. With an appendix by David Mumford. Ergebnisse der Mathematik und ihrer Grenzgebiete (3), 22. Springer-Verlag, Berlin, 1990. 

\bibitem{GIT} J. Fogarty, F. Kirwan, D. Mumford: Geometric invariant theory, Ergebnisse der Mathematik und ihrer Grenzgebiete. 2. Folge, 34. Springer-Verlag Berlin Heidelberg, 1994.

\bibitem{vG98} B. van Geemen: \emph{The Schottky problem and second order theta functions}. Workshop on Abelian Varieties and Theta Functions (Spanish) (Morelia, 1996), 41--84,
Aportaciones Mat. Investig., 13, Soc. Mat. Mexicana, México, 1998.

\bibitem{vG} B. van Geemen: \emph{Some equations for the universal Kummer variety}. Trans. Amer. Math. Soc. 368 (2016), no. 1, 209--225. 

\bibitem{vdG} G. van der Geer: \emph{On the geometry of a Siegel modular threefold}. Math. Ann. 206 (1982) no.3, 317--350.

\bibitem{GD} M.R. Gonzalez-Dorrego: \emph{$(16,6)$ configurations and geometry of Kummer surfaces in $\PP^3$}. Mem. Amer. Math. Soc. 107 (1994), no. 512.

\bibitem{M2}
D. Grayson, M. Stillman, \textit{Macaulay2, a software system for research in algebraic geometry}. 

\bibitem{G}
S. Grushevsky, K. Hulek: \emph{Geometry of theta divisors--a survey}. In: A celebration of algebraic geometry. 361–390, Clay Math. Proc., 18, Amer. Math. Soc., Providence, RI, 2013.

\bibitem{GH2011} S. Grushevsky, K. Hulek: \emph{Principally Polarized Semi-abelic Varieties of Small
Torus Rank, and the Andreotti-Mayer Loci}. Pure and Applied Mathematics Quarterly, Volume 7, Number 4, 1309--1360, 2011.


\bibitem{ig1} 
J.-I. Igusa: 
\textit{Theta functions}. 
Grundlehren der Mathematischen Wissenschaften, 194. 
Springer-Verlag, New York-Heidelberg, 1972.

\bibitem{Meq1}
D. Mumford: \emph{On the Equations Defining Abelian Varieties I}. Invent. math. 1, 287--354 (1966).

\bibitem{Mum70} D. Mumford: \emph{Abelian Varieties}. Tata Institut of Fundamental
Research Studies in Mathematics, No. 5. Published for the Tata
Institute of Fundamental Research, Bombay, 1970.

\bibitem{Mum83} D. Mumford: \emph{On the Kodaira dimension of the Siegel modular variety}. Lect. Notes in Math., Vol. 997 (1983),
348-375, Springer.



\bibitem{RSSS13} Q. Ren, S. V. Sam, G. Schrader, B. Sturmfels: \emph{The Universal Kummer Threefold}. Exp. Math. 22 (2013), no. 3, 327–362. 

\bibitem{roma05}
M. Romagny: \emph{Group actions on stacks and applications}. Michigan Math. J., 53(1):209--236, 2005.

\bibitem{SB16} N. I. Shepherd-Barron: \emph{An exceptional locus in the perfect compactification of $ A_g$}.  arXiv preprint arXiv:1604.05954, 2016

\bibitem{S-B} N. I. Shepherd-Barron: \emph{Perfect forms and the moduli space of abelian varieties}. Invent. Math., 163(1):25--45, 2006.



\end{thebibliography}
\end{document}